
\input amstex
\magnification=1200
\font\cyr=wncyr10
\font\cyi=wncyi10
\font\cyb=wncyb10
\font\cyre=wncyr8
\font\cyie=wncyi8
\documentstyle{amsppt}
\NoRunningHeads
\NoBlackBoxes
\hyphenation{izo-to-pi-ches-kaya izo-to-pi-ches-kih ope-ra-tor}

\define\End{\operatorname{End}}
\define\Mat{\operatorname{Mat}}
\define\gl{\operatorname{gl}}
\define\osp{\operatorname{osp}}
\define\osq{\operatorname{osq}}
\define\q{\operatorname{q}}
\define\Vect{\operatorname{Vect}}
\define\sLtwo{\operatorname{sl}(2,\Bbb C)}
\define\Ind{\operatorname{Ind}}
\define\Map{\operatorname{Map}}
\document\qquad\qquad\qquad\qquad\qquad\qquad\qquad\qquad\qquad
$\boxed{\boxed{\aligned
&\text{\eightpoint"Thalassa Aitheria" Reports}\\
&\text{\eightpoint RCMPI-95/06}\endaligned}}$\newline
\ \newline
\ \newline
\ \newline
\topmatter
\title\cyb Izotopicheskie pary i ih predstavleniya. {\rm II}. Obshchii0
supersluchai0.
\endtitle
\author\cyr D.V.Yurp1ev
\endauthor
\address\cyre CMFI "Talassa E1teriya"
\endaddress
\email denis\@juriev.msk.ru\endemail
\abstract\nofrills\eightpoint\cyre\newline
Izotopicheskie pary, algebraicheskie obp2ekty, posredestvom kotoryh
pre\-d\-stav\-lya\-et\-sya udobnym opisyvatp1 nekotorye formy negamilp1tonova
vza\-imo\-dee0\-st\-viya (magnitnogo tipa) gamilp1tonovyh sistem na kvantovom
urov\-ne, ras\-smat\-ri\-va\-yut\-sya v maksimalp1noe0 obshchnosti v ramkah
formalizma vek\-tor\-noe0 superalgebry.
\endabstract
\endtopmatter
\cyr

\

\

V pervoe0 chasti raboty, posvyashchennoe0 izotopicheskim param i ih
predstavleniyam, rassmatrivalisp1 chisto chetnye izotopicheskie pary [1].
Od\-na\-ko, kak s teoreticheskoe0, tak i s prikladnoe0 tochek zreniya imeet
smysl issledovatp1 naryadu s chisto chetnymi obp2ektami ih nechetnye i
superanalogi [2]. V dannoe0 vtoroe0 chasti izotopicheskie pary
izu\-cha\-yut\-sya vo vsee0 obshchnosti.

\subhead 1. \cyb Izotopicheskie pary: opredeleniya i primery\endsubhead
\cyr Oboznachim
$$\aligned A_{XYZ}=&(-1)^{p(X)p(Y)+p(Y)p(Z)+p(Z)p(X)},\\
B_{XYZW}=&(-1)^{p(X)p(Y)+p(Y)p(Z)+p(Z)p(W)+p(W)p(X)},\endaligned$$
gde $p(X)=0, 1$ --- chetnostp1 e1lementa $X$.

\definition{\cyb Opredelenie 1} \cyr Para $(V_1,V_2)$ linee0nyh
superprostranstv $V_i=V_i^0\oplus V_i^1$ ($i=1,2$) nazyvaet\-sya {\cyi
izotopicheskoe0 paroe0} ({\it isotopic pair\/}), esli i tolp1ko esli
opredeleny dva (chetnyh) otobrazheniya
$m_1:V_2\otimes V_1\otimes V_1\mapsto V_1$ i
$m_2:V_1\otimes V_2\otimes V_2\mapsto V_2$ takie, chto operatsii $(X,Y)\mapsto
[X,Y]_U=m_i(U,X,Y)=A_{XUY}m_i(U,Y,X)$ ($i=1,2$; $X,Y\in V_1$,
$U\in V_2$ ili $X,Y\in V_2$, $U\in V_1$), nazyvaemye izokommutatorami
(s izotopicheskimi e1lementami $U$), udovletvoryayut tozhdestvam
$$\aligned &A_{VYZ}[[X,Y]_U,Z]_V+A_{UZV}[X,[Z,Y]_U]_V+A_{XZU}[[Z,X]_U,Y]_V+\\
&B_{VZYU}[[X,Y]_V,Z]_U+[X,[Z,Y]_V]_U+B_{XZUV}[[Z,X]_V,Y]_U=0
\endaligned$$
(analogam tozhdestva Yakobi) i soglasovany mezhdu soboe0:
$$\aligned
&[X,Y]_{[U,V]_Z}=\tfrac12(A_{VYZ}[[X,Y]_U,Z]_V+[X,[Z,Y]_V]_U+
B_{XZUV}[[Z,X]_V,Y]_U)-\\
&\tfrac12(B_{VZYU}[[X,Y]_V,Z]_U+A_{UZV}[X,[Z,Y]_U]_V+
A_{XZU}[[Z,X]_U,Y]_V)\endaligned$$
gde $X,Y,Z\in V_1$, $U,V\in V_2$ ili $X,Y,Z\in V_2$, $U,V\in V_1$.
Para $(V_1,V_2)$ linee0nyh superprostranstv
$V_i=V_i^0\oplus V_i^1$ ($i=1,2$) nazyvaet\-sya {\cyi
su\-per\-e0or\-da\-no\-voe0 paroe0} ({\it super--Jordan pair\/}), esli i
tolp1ko esli opredeleny dva (chetnyh) otobrazheniya
$m_1:V_2\otimes V_1\otimes V_1\mapsto V_1$ i
$m_2:V_1\otimes V_2\otimes V_2\mapsto V_2$ takie, chto operatsii $(X,Y)\mapsto
X\underset{U}\to\circ Y=m_i(U,X,Y)=A_{XUY}m_i(U,Y,X)$ ($i=1,2$; $X,Y\in V_1$,
$U\in V_2$ ili $X,Y\in V_2$, $U\in V_1$) soglasovany mezhdu soboe0
sleduyushchim obrazom:
$$\align
X\underset{U\underset{Z}\to\circ V}\to\circ Y=
&X\underset{U}\to\circ(Z\underset{V}\to\circ Y)-
A_{VYZ}(X\underset{U}\to\circ Y)\underset{V}\to\circ Z-
B_{XZUV}(Z\underset{V}\to\circ X)\underset{U}\to\circ Y
\endalign$$
gde $X,Y,Z\in V_1$, $U,V\in V_2$ ili $X,Y,Z\in V_2$, $U,V\in V_1$.
Superalgebra Li $\frak g$ takaya, chto para superprostranstv
$(\frak g,\frak g)$ nadelena $\frak g$--e1kvivariantnoe0 strukturoe0
izotopicheskoe0 pary s izokommutatorami, soglasovannymi s is\-hodnymi
superskobkami Li v $\frak g$, nazyvaet\-sya {\cyi magnitnoe0 superalgebroe0
Li\/} ({\it magnetic Lie algebra\/)}.
\enddefinition

\cyr
Obsudim e1to opredelenie.

Vo-pervyh, otmetim, chto dannoe opredelenie izotopicheskih par
yav\-lya\-et\-sya bolee obshchim po otnosheniyu k opredeleniyu rabot [1,3], v
kotoryh rassmatrivalsya chisto chetnye0 sluchae0.

Vo-vtoryh, opredelenie izotopicheskoe0 pary mozhet rassmatrivatp1sya kak
rezulp1tat aksiomatizatsii sleduyushchee0 konstruktsii: pustp1 $\Cal A$ ---
assotsiativnaya superalgebra, napr. matrichnaya, $V_1$, $V_2$ --- dva linee0nyh
podprostranstva v nee0 takie, chto $V_1$, $V_2$ zamknuty otnositelp1no
izokommutatorov $(X,Y)\mapsto [X,Y]_U=XUY-A_{XUY}YUX$ s izotopicheskimi
e1lementami $U$ iz $V_2$, $V_1$, sootvet\-stvenno. Opredelenie
supere0ordanovoe0 pary po\-lu\-cha\-et\-sya analogichnym obrazom, esli vmesto
izokommutatorov rassmatrivatp1 operatsii $X\underset{U}\to{\circ}
Y=XUY+A_{XUY}YUX$. Ukazannaya konstruktsiya proyasnyaet nazvanie
"izotopicheskaya para": v samom dele, izokommutator $[X,Y]_U$ udovletvoryaet
sootnosheniyam $Q[X,Y]_Q=[QX,QY]$, $[X,Y]_QQ=[XQ,YQ]$, inymi slovami, izotopen
obychnomu (super)kommutatoru. Izotopii v e1tom sluchae ne obyazatelp1no
nevyrozhdeny (t.e. $Q^{-1}$ ne vsegda sushchestvuet), no linee0ny, i
zadayut\-sya pri pomoshchi levogo ili pravogo umnozheniya v algebre $A$:
$X\mapsto QX$ ili $X\mapsto XQ$. Otmetim, chto obshchie izotopii algebr imeyut
bolee slozhnye0 vid, chem ispolp1zuemye v dannoe0 konstruktsii [4]. Naprimer,
imeet smysl rassmatrivatp1 takzhe nelinee0nye kvadratichnye izotopii vida
$X\mapsto QXQ$ (sr.[5]), regulyarno voznikayushchie v teorii e0ordanovyh
algebr, vypuklyh konusov i simmetricheskih prostranstv [6].

V-tretp1ih, esli $\frak g$ -- magnitnaya superalgebra Li, to
$(\frak g\oplus\Bbbk,\frak g\oplus\Bbbk)$ -- izotopicheskaya para. Naoborot,
lyubaya izotopicheskaya para $(V_1,V_2)$ takaya, chto $\dim V_1=\dim V_2=n|m$,
yavlyaet\-sya magnitnoe0 superalgebroe0 Li s abelevoe0 superalgebroe0 Li
$\frak g=\Bbbk^{n|m}$.

V-chetvertyh, otmetim, chto chetnye supere0ordanovy pary yavlyayut\-sya
e0o\-r\-da\-no\-vy\-mi [7], a nechetnye supere0ordanovy pary ---
an\-ti\-e0or\-da\-no\-vy\-mi [8].

V-pyatyh, esli harakteristika osnovnogo polya bolp1she dvuh, to izotopicheskie
pary nadelyayut\-sya strukturoe0 supere0ordanovyh par pri izmenenii
supergraduirovki na protivopolozhnuyu, i {\rm vice versa}. Takim obrazom, v
e1tom sluchae chisto chetnye izotopicheskie pary otozhdestvlyayut\-sya s
antie0ordanovymi parami. Odnako, antie0ordanovy pary nad polem harakteristiki
dva predostavlyayut primer takih par, ne yavlyayushchihsya izotopicheskimi
(sr. zamechaniya v [1,3]).

Razberem primery izotopicheskih par.

Vo-pervyh, esli $\frak g$ --- superalgebra Li, to para $(\frak g,\Bbbk)$
nadelyaet\-sya estestvennoe0 strukturoe0 izotopicheskoe0 pary. Bolee togo,
nekotorye (li\-nee0\-nye) puchki superalgebr Li $\frak g_{\lambda}$ nad polem
$\Bbbk$ ($\lambda\in\Bbbk^n$) mogut rassmatrivatp1sya kak izotopicheskie pary
$(V_1,V_2)$, gde $V_1\simeq\frak g_{\lambda}$, $V_2\simeq\Bbbk^n$, pri
vypolnenii nekotoryh uslovie0 soglasovaniya na superkommutatory
($\lozenge$--zamknutostp1). Chisto chetnye0 $\frak g$--e1kvivariantnye0
sluchae0 (lievy $\frak g$--puchki) izuchalsya v rabotah [9,par.3; 3,app.A],

Vo-vtoryh, mnogochislennye primery chetnyh izotopicheskih
(an\-ti\-e0or\-da\-no\-vyh) par rassmatrivalisp1 v rabotah [1,3]. Nekotorye iz
nih ot\-ve\-cha\-yut magnitnym algebram Li (t.e. chisto chetnym magnitnym
superalgebram Li).

V-tretp1ih, imeet mesto sleduyushchaya teorema.

\proclaim{\cyb Teorema 1A} \cyi
Sushchestvuet pyatp1 beskonechnyh serie0 izotopicheskih par:

(1) izotopicheskaya para $\gl(n,m)$, obrazovannaya matritsami
$\left(\matrix A & B \\ C & D\endmatrix\right)$ i\linebreak
$\left(\matrix X & Y \\ Z & W\endmatrix\right)$ iz $\Mat(n|m)$;

(2,3) izotopicheskie pary $\osp^{\pm}(n,m)$, yavlyayushchiesya podparami
izo\-to\-pi\-ches\-koe0 pary $\gl(n,m)$, vydelyaemymi usloviyami:
$A^t=-A$, $D^t=D$, $B^t=\pm C$, $X^t=X$, $W^t=-W$, $Y^t=\pm Z$;

(4) izotopicheskaya para $\q(n)$, yavlyayushchayasya podparoe0 izotopicheskoe0
pary $\gl(n,n)$, vydelyaemoe0 usloviyami: $A=D$, $B=C$, $X=W$, $Y=Z$;

(5) izotopicheskaya para $\osq(n)$, yavlyayushchayasya podparoe0
izotopicheskoe0 pary $\q(n)$, vydelyaemoe0 usloviyami: $A^t=A$, $B^t=-B$,
$X^t=-X$, $W^t=W$.
\endproclaim

\cyr
Konechno zhe, utverzhdenie teoremy ne ogranichivaet klass izotopicheskih par
pyatp1yu perechislennymi seriyami. Naprimer, dlya lyubyh dvuh e1lementov $A$ i
$B$ izotopicheskoe0 pary $(V_1,V_2)$ ($A\in V_1$, $B\in V_2$) prostranstva
$V_1^{\dag}=\{X\in V_1: [A,X]_B=0\}$ i $V_2^{\dag}=\{Y\in V_2: [B,Y]_X=0\}$
obrazuyut izotopicheskuyu podparu $(V_1^{\dag},V_2^{\dag})$ pary $(V_1,V_2)$
(dannaya konstruktsiya predstavlyaet interes v kontekste raboty [10]). Chisto
chetnye izotopicheskie pary ukazannyh serie0 (a takzhe drugie primery)
rassmatrivalisp1 v rabote [1].

V-chetvertyh, ime\-et\-sya ryad netrivialp1nyh beskonechnomernyh izotopicheskih
par.

\proclaim{\cyb Teorema 1B} \cyi Pary $(W(n|m),O(n|m))$,
$(K(2n+1|m),O(2n+1|m))$,
$(M(n),O(n|\mathbreak n+1))$, $(H_Ó(2n|m),O(2n|m))$,
$(Le_Ó(n),O(n|n))$ yavlyayut\-sya izotopicheskimi pa\-ra\-mi, gde $O(n|m)$ --
prostranstvo formalp1nyh stepennyh ryadov ot $n$ chetnyh i $m$ nechetnyh
peremennyh, $W(n|m)$, $K(2n+1|m)$, $M(n)$, $H_Ó(2n|m)$, $Le_Ó(n)$ ---
superalgebra Li formalp1nyh vektornyh polee0, kontaktnaya i nechetnaya
kontaktnaya superalgebry Li, superalgebra konformno--gamilp1tonovyh vektornyh
polee0 i ee nechetnye0 analog [11], sootvetstvenno.

Pary $(\Vect(M),\Cal O(M))$ É $(\Vect(M)\oplus\Cal O(M),\Cal O(M))$,
gde $M$ --- proizvolp1noe supermnogoobrazie, $\Cal O(M)$ i $\Vect(M)$ ---
prostranstva funktsie0 i vektornyh polee0 na nem, so\-ot\-ve\-t\-s\-t\-ven\-no,
yavlyayut\-sya izotopicheskimi parami, pri\linebreak e1tom mozhno
ogranichitp1sya kontaktnymi ili kon\-for\-m\-no--ga\-milp1\-to\-no\-vy\-mi
vektornymi polyami (ispolp1zuya v sluchae pervoe0 pary otobrazhenie
prostranstva so\-ot\-vet\-s\-t\-vu\-yu\-shchih vektornyh polee0 na prostranstvo
dif\-fe\-ren\-tsi\-alp1\-nyh operatorov pervogo poryadka, annuliruyushchih
so\-ot\-vet\-s\-t\-vu\-yu\-shchuyu\linebreak formu). V chastnosti,
izotopicheskimi parami yavlyayut\-sya $(\Cal K(n),\Cal O(S^{1|n}))$ i
$(\Cal K^+(n),\Cal O(S^{1|n}_+))$, gde $\Cal K(n)$ i $\Cal K^+(n)$ ($n=1,2,3$)
--- superalgebry Neve--Shvartsa i Ramona (pri $n=1$) i ih vysshie analogi,
faktorizovannye po tsentru, a $\Cal O(S^{1|n})$ i $\Cal O(S^{1|n}_+)$ ---
prostranstva funktsie0 na superokruzhnostyah $S^{1|n}$ i podkruchennyh na
rassloenie Mebiusa superokruzhnostyah $S^{1|n}_+$ (smotri [11]).
\endproclaim

\cyr
Bylo by interesno opisatp1 tsentralp1nye rasshireniya
bes\-ko\-ne\-chno\-mer\-nyh izotopicheskih par teoremy 1B, a takzhe
"izotopicheskih par tokov" vida $(\Map(\Bbb S^1,V_1),\Map(\Bbb S^1,V_2))$
($(V_1,V_2)$ -- izotopicheskaya para), ih mno\-go\-mer\-nyh i
superobobshchenie0 ($\Bbb S^1$ zamenyaet\-sya na proizvolp1noe
supermnogoobrazie $\Cal M$).

Kak otmechalosp1 vyshe, primery magnitnyh algebr Li soderzhat\-sya v rabote
[1]. Ukazhem takzhe, chto poluprostaya algebra Li $\frak g$ yavlyaet\-sya
magnitnoe0 algebroe0 Li s izokommutatorami $[X,Y]_U=\pm(X,U)Y\mp(U,Y)X$, gde
$(\cdot,\cdot)$ --- bilinee0naya forma Killinga. Magnitnoe0 algebre Li $\sLtwo$
otvechaet izokvaternionnaya izotopicheskaya para (sm. [1]).

Krome togo, privedem prostuyu konstruktsiyu, dopuskayushchuyu neposredstvennoe
superobobshchenie. A imenno, esli $\frak g$ --- algebra Li s invariantnoe0
simmetricheskoe0 bilinee0noe0 formoe0 $\eta_{\alpha\beta}$ i strukturnymi
konstantami $c_{\alpha\beta}^{\gamma}$, $c_{\alpha\beta\gamma}=
1c_{\alpha\beta}^{\delta}\eta_{\gamma\delta}$, to $(\frak g, S^2(\frak g))$ ---
izotopicheskaya para s izokommutatorami:
$$\aligned
[e_{\alpha},e_{\beta}]_{m_{\gamma\delta}}=&(\eta_{\alpha\gamma}e^{\rho}_{\beta
\delta}-\eta_{\alpha\gamma}c^{\rho}_{\alpha\delta})e_{\rho},\\
[m_{\alpha\beta},m_{\gamma,\delta}]_{e_{\zeta}}=&c_{\beta\gamma\zeta}
m_{\alpha\delta}+c_{\alpha\gamma\zeta}m_{\beta\delta}+c_{\beta\delta\zeta}
m_{\alpha\gamma}+c_{\alpha\delta\zeta}m_{\beta\gamma}.
\endaligned
$$
gde $\{e_{\alpha}\}$, $\{m_{\alpha\beta}=m_{\beta\alpha}\}$ --- bazisy v
$\frak g$ i $S^2(\frak g)$, sootvetstvenno. Pri e1tom, $(\frak g,
S^2(\frak g)/S^2(\frak g)^{\frak g})$ takzhe yavlyaet\-sya izotopicheskoe0
paroe0.

\subhead 2. \cyb Predstavleniya izotopicheskih par: opredeleniya i primery
\endsubhead

\definition{\cyb Opredelenie 2} \cyr {\cyi Predstavleniem izotopicheskoe0
pary $(V_1,V_2)$ v linee0nom superprostranstve $H$\/} nazyvaet\-sya para
$(T_1,T_2)$ (chetnyh) otobrazhenie0 $T_i:V_i\mapsto\End(H)$ takih, chto
$$\align
T_1([X,Y]_U)=&T_1(X)T_2(U)T_1(Y)-A_{XUY}T_1(Y)T_2(U)T_1(X),\\
T_2([U,V]_X)=&T_2(U)T_1(X)T_2(V)-A_{UXV}T_2(U)T_1(X)T_2(V),
\endalign$$
gde $X,Y\in V_1$, $U,V\in V_2$. Predstavlenie izotopicheskoe0 pary $(V_1,V_2)$
v linee0nom superprostranstve $H$ nazyvaet\-sya {\cyi rasshcheplennym\/} ({\it
split\/}), esli i tolp1ko esli $H=H_1\oplus H_2$ i
$$\left\{\aligned
(\forall X\in V_1) \left. T_1(X)\right|_{H_2}=0,\ T_1(X):H_1\mapsto H_2,\\
(\forall U\in V_2) \left. T_2(U)\right|_{H_1}=0,\ T_2(U):H_2\mapsto H_1.
\endaligned\right.$$
\enddefinition

\define\vac{{\left|0\right>}}

\remark{\cyi Primer 1. Rasshcheplennye predstavleniya so starshim vesom
(\it Highest weight split representations)}

\cyr
Esli izotopicheskaya para $(V_1,V_2)$ $\Bbb Z$--graduirovana (t.e.
$V_1=\bigoplus_{i\in\Bbb Z}
V_{1,i}$, $V_2=\bigoplus_{i\in\Bbb Z} V_{2,i}$,
$[V_{1,i},V_{1,j}]_{V_{2,k}}\subseteq V_{1,i+j+k}$,
$[V_{2,i},V_{2,j}]_{V_{1,k}}\subseteq V_{2,i+j+k}$)
i izotopicheskaya para $(V_{1,0},V_{2,0})$ trivialp1na, to mozhno
rassmatrivatp1 rasshcheplennye predstavleniya pary $(V_1,V_2)$ so starshim
vesom. Prostranstva $H_1$ i $H_2$ porozhdayut\-sya starshimi (vakuumnymi)
vektorami $\vac_1\in H_1$ i $\vac_2\in H_2$ takimi, chto
$(\forall X\in V_{1,0}) X\vac_1=\chi_1(X)\vac_2$,
$(\forall U\in V_{2,0}) X\vac_2=\chi_2(U)\vac_1$,
gde $\chi_i\in V^*_{i,0}$,
$(\forall X\in V_{1,-i}) X\vac_1=0$,
$(\forall U\in V_{2,-i}) X\vac_2=0$ ($i>0$).
Ostalp1nye bazisnye vektory v $H_1$ i $H_2$ imeyut vid
$U_jX_j\ldots U_1X_1\vac_1$ ili $U_{j+1}X_jU_j\ldots\mathbreak X_1U_1\vac_2$ i
$X_jU_j\ldots X_1U_1\vac_2$ ili $X_{j+1}U_jX_j\ldots U_1X_1\vac_1$,
sootvet\-stvenno, gde $X_1,\ldots X_j,X_{j+1}\in\bigoplus_{i>0}V_{1,i}$,
$U_1,\ldots U_j,U_{j+1}\in\bigoplus_{i>0}V_{2,i}$.
\endremark

\remark{\cyi Primer 2. Indutsirovannye rasshcheplennye predstavleniya}

\cyr
Konstruktsiya primera 1 neposredstvenno obobshchaet\-sya, predostavlyaya
ana\-lo\-gi indutsirovannyh predstavlenie0.

\define\vc{{\left|v\right>}}
A imenno, pustp1 $(V_1,V_2)$ --- izotopicheskaya para, $(V_1^\circ,V_2^\circ)$
--- ee izotopicheskaya podpara i $T^\circ=(T_1^\circ,T_2^\circ)$ ---
predstavlenie $(V^\circ_1,V^\circ_2)$ v prostranstve $H^\circ=H^\circ_1\oplus
H^\circ_2$. Togda mozhno rassmotretp1 indutsirovannoe rasshcheplennoe
predstavlenie
$T=(T_1,T_2)=\operatorname{Ind}^{(V_1,V_2)}_{(V_1^\circ,V_2^\circ)}(T_1^\circ,
T_2^\circ)$ pary $(V_1,V_2)$ v prostranstve $H=H_1\oplus H_2$, gde $H_1$ ---
porozhdeno $U_jX_j\ldots U_1X_1\vc_1$ i $U_{j+1}X_jU_j\ldots X_1U_1\vc_2$, v to
vremya kak $H_2$ porozhdeno $X_jU_j\ldots X_1U_1\vc_2$ i
$X_{j+1}U_jX_j\ldots U_1X_1\vc_1$\linebreak ($X_1,\ldots X_j,X_{j+1}\in V_1$;
$U_1,\ldots U_j,U_{j+1}\in V_2$; $\vc_1\in H^\circ_1$, $\vc_2\in
H^\circ_2$).
\endremark

\cyr
Predstavleniya so starshim vesom i indutsirovannye rasshcheplennye
pre\-d\-s\-tav\-le\-niya mogut bytp1 postroeny dlya izotopicheskih par pyati
serie0 teoremy 1A. Otmetim takzhe, chto mozhno rassmatrivatp1 i konstruktsiyu
geometricheskih predstavlenie0 po analogii s [1].

K sozhaleniyu, avtoru neizvestno, chto sleduet s\-chitatp1 (ko)gomologiyami
izotopicheskih par, no e1ti obp2ekty, nadlezhashchim obrazom opredelennye,
nesomnenno dolzhny predstavlyatp1 interes (sr.[12]).

Predstavleniya izotopicheskih par $(\frak g,\Bbbk)$, assotsiirovannye s
algebrami Li $\frak g$, rassmatrivalisp1 v rabote [3]. Po kazhdomu
predstavleniyu $T=(T_1,T_2)$ izotopicheskoe0 pary $(\frak g,\Bbbk)$ mozhet
bytp1 postroeno predstavlenie $T_0$ superalgebry Li $\frak g$:
$T_0(X)=T_2(1)T_1(X)$, $1\in\Bbbk$, $X\in\frak g$. Naoborot, esli $T_0$ ---
nekotoroe predstavlenie superalgebry Li $\frak g$ i $Q$ --- nevyrozhdennye0
operator v prostranstve predstavleniya, to $(T_1,T_2)$, gde
$T_1(X)=Q^{-1}T_0(X)$, $T_2(1)=Q$ ($1\in\Bbbk$, $X\in\frak g$), yavlyaet\-sya
predstavleniem izotopicheskoe0 pary $(\frak g,\Bbbk)$. V sluchae vyrozhdennogo
$Q$ mogut estestvennym obrazom voznikatp1 i bolee slozhnye nelievskie
algebraicheskie obp2ekty, takie kak, naprimer, "kvadratichnye" algebry
Zhelobenko ($AZ_n$, $BZ_n$, $CZ_n$, $DZ_n$), Mikelp1sona ($S(\frak g,\frak k)$,
$Z(\frak g,\frak k)$) i ih obobshcheniya ($Z(A,\frak k)$, gde $A$, naprimer,
--- kontragredientnaya assotsiativnaya algebra Zhelobenko) [13], esli $Q$ ---
t.n. al\-geb\-rai\-ches\-kie0 e1kstremalp1nye0 proektor
Asherovoe0--Smirnova--Tolstogo [14,13].

\subhead 3. \cyb Nekotorye algebraicheskie obp2ekty, svyazannye s
izotopicheskimi i supere0ordanovymi parami i magnitnymi superalgebrami
\endsubhead
\cyr
Horosho izvestno, po e0ordanovoe0 pare $(V_1,V_2)$ v prostranstve $V=V_1\oplus
V_2$ stroit\-sya polyarizovannaya troe0naya sistema Li [7], a po
anti-e0ordanovoe0 pare (i, v chastnosti, po chisto chetnoe0 izotopicheskoe0
pare) $(V_1,V_2)$ v prostranstve $V=V_1\oplus V_2$ stroit\-sya polyarizovannaya
antilieva troe0naya sistema [8].

V svoyu ocheredp1 ukazannym troe0nym sistemam sopostavlyaet\-sya nekotorye0
klass algebr i superalgebr Li.

Analogichnaya situatsiya, okazyvaet\-sya, imeet mesto i v obshchem
supersluchae.

\proclaim{\cyb Teorema 2A} \cyi Pustp1 $(V_1,V_2)$ -- supere0ordanova para,
togda prostranstvo $V=V_1\oplus V_2$ nadelyaet\-sya strukturoe0
polyarizovannoe0 superlievoe0 troe0noe0 sistemy (opredelenie smotri v [15,
str.2786]).
\endproclaim

\cyr
Konstruktsiya polyarizovannoe0 superlievoe0 troe0noe0 sistemy v tochnosti
povtoryaet analogichnye konstruktsii v chisto chetnom i nechetnom sluchayah.
Tem samym, dokazatelp1stvo teoremy svodit\-sya k prostoe0 identifikatsii
opredelyayushchih sootnoshenie0 v oboih klassah algebraicheskih sistem.

Takim obrazom, po kazhdoe0 izotopicheskoe0 pare mozhet bytp1 postroena
polyarizovannaya superlieva troe0naya sistema. Primery dlya chisto chetnogo
sluchaya rassmatrivalisp1 vesp1ma podrobno v [3].

Budem nazyvatp1 $\Bbb Z_2$--graduirovannuyu superalgebru Li $\frak g=\frak g_0
\oplus\frak g_1$ po\-lya\-ri\-zo\-van\-noe0, esli linee0noe prostranstvo
$\frak g_1$ dopuskaet razlozhenie $\frak g_1=\frak g_1^+\oplus\frak g_1^-$
takoe, chto (1) $\frak g_1^{\pm}$ -- $\frak g_0$--podmoduli modulya
$\frak g_1$, (2) $[\frak g_1^+,\frak g_1^+]=[\frak g_1^-,\frak g_1^-]=0$.
Tselesoobrazno rassmatrivatp1 dopolnitelp1nuyu $\Bbb Z_2$--graduirovku kak
"podkruchivanie" is\-hodnoe0 supergraduirovki, t.e. snabzhatp1 $\frak g_1$
pro\-ti\-vo\-po\-lozh\-noe0 chetnostp1yu.

\proclaim{\cyb Teorema 2B} \cyi Kazhdoe0 izotopicheskoe0 pare $(V_1,V_2)$
sopostavlyaet\-sya polyarizovannaya $\Bbb Z_2$--graduirovannaya superalgebra Li
$\frak g=\frak g_0\oplus\frak g_1=\frak g_0\oplus(\frak g_1^+\oplus\frak
g_1^-)$
takaya, chto $\frak g_1^+=V_1$ É $\frak g_1^-=V_2$.
\endproclaim

\cyr
Konstruktsiya kopiruet chisto chetnye0 sluchae0, podrobno
pro\-kom\-men\-ti\-ro\-van\-nye0 i proillyustrirovannye0 na primerah v rabote
[3]. Podcherknem, chto esli v pervoe0 chasti teoremy 2 estestvenno
rassmatrivatp1 su\-per\-e0or\-da\-no\-vy pary, to vo vtoroe0 ---
izotopicheskie.

Magnitnye poluprostye superalgebry Li $\frak g$ odnoznachno opredelyayut\-sya
nulevoe0 komponentoe0 $\hat{\frak g}_0$ sootvet\-stvuyushchee0 ee0 po teoreme
2B superalgebry Li $\hat{\frak g}$, monomorfizmom $\frak g\oplus\frak g
\mapsto\hat{\frak g}_0$ i razlozheniem $\hat{\frak g}_0=
(\frak g\oplus\frak g)\oplus\frak w$, gde $\frak w$ -- podalgebra v
$\hat{\frak g_0}$, porozhdennaya izokommutatorami. Podobnye razlozheniya
igrayut vazhnuyu rolp1 v formalizme klassicheskoe0 $r$--matritsy [16].

\proclaim{\cyb Teorema 3A} \cyi V prostranstve rasshcheplennogo predstavleniya
izo\-to\-pi\-ches\-koe0 pary $(V_1,V_2)$ realizuet\-sya predstavlenie
sootvet\-stvuyushchee0 polyarizovannoe0 $\Bbb Z_2$--graduirovannoe0
superalgebry Li.
\endproclaim

\cyr
Chisto chetnye0 analog e1toe0 teoremy priveden v [3].

\proclaim{\cyb Teorema 3B} \cyi Pustp1 $(V_1,V_2)$ --- nekotoraya
izotopicheskaya para s re\-duk\-tiv\-noe0 so\-ot\-vet\-s\-t\-vu\-yu\-shchee0
superalgebroe0 Li $\frak g$, $\frak p$ --- parabolicheskaya podalgebra v
$\frak g$ (sr. [17]) takaya, chto $\dim(\frak g_1/(\frak g_1\cap\frak
p)=(0|1)$,
$\chi$ --- harakter $\frak p$, togda v prostranstve predstavleniya
$\Ind_{\frak p}^{\frak g}\chi$ superalgebry Li $\frak g$, indutsirovannogo s
haraktera $\chi$ podalgebry $\frak p$, realizuet\-sya rasshcheplennoe
predstavlenie izo\-to\-pi\-ches\-koe0 pary $(V_1,V_2)$.
\endproclaim

\cyr
V sluchae izokvaternionnoe0 pary (otvechayushchee0 magnitnoe0 algebre Li
$\sLtwo$) konechnomernye podpredstavleniya indutsirovannyh predstavlenie0
opisyvayut vnutrennyuyu sverhtonkuyu magnitnuyu strukturu svyazannyh
so\-sto\-ya\-nie0 dvuh chastits so spinom. V chastnosti, chetyrehmernoe
rasshcheplennoe predstavlenie so starshim vesom $(1/2,1/2)$ (fundamentalp1noe
predstavlenie izokvaternionnoe0 pary) sootvet\-stvuet pare svyazannyh chastits
spina $1/2$ (sm. napr. [18]). Krome togo, te zhe dannye harakterizuyut
vnutrennyuyu sverhtonkuyu strukturu skalyarno--tenzornogo
vzaimodee0stviya\linebreak nuklonov v nuklon--nuklonnyh parah, vzaimodee0stviya
nuklon--nuklonnyh par, pozitroniev ili kuperovskih par mezhdu soboe0, a takzhe
nekotoryh form negamilp1tonova vzaimodee0stviya kvantovyh vihree0 v
sverhprovodnikah ili kvantovyh zhidkostyah. V silu e1togo predstavlyaet
interes obshchaya zadacha izucheniya spektra ogranichenie0 neprivodimyh
predstavlenie0 izotopicheskih par $(\frak g\oplus\Bbbk,\frak g\oplus\Bbbk)$
($\frak g$ --- magnitnaya superalgebra Li) i sootvet\-stvuyushchih
$\Bbb Z_2$--graduirovannyh polyarizovannyh superalgebr Li $\hat{\frak g}$
na superalgebru Li $\frak g\oplus\frak g\subseteq\hat{\frak g}_0$ i pravil
superotbora (sm. napr. [19]). Dannaya zadacha imeet smysl i dlya
beskonechnomernyh izotopicheskih par (pe\-re\-chis\-len\-nyh v teoreme 2B i
"izotopicheskih par tokov") vvidu vozmozhnyh prilozhenie0 k teorii (super)strun
(fenomenologiya chastits i spontannoe narushenie simmetrii) i strunnopodobnyh
kvantovyh obp2ektov.

\subhead 4. \cyb Nekotorye zamechaniya\endsubhead
\cyr
Izotopicheskie pary (a takzhe magnitnye algebry Li) yavlyayut\-sya
prostee0shimi algebraicheskimi obp2ektami, opisyvayushchimi nepotentsialp1noe
vzaimodee0stvie gamilp1tonovyh sistem. Bolee obshchie struktury --- $I$--pary
[5], v kotoryh vzaimodee0stvie ne obyazatelp1no linee0no, a takzhe obshchih
nelinee0nyh kvantovyh (ili klassicheskih) skobok (sr. napr. [20,21]),
zavisyashchih ot sostoyanie0 vozdee0stvuyushchih sistem. Podobnye obp2ekty
mogut vstrechatp1sya, naprimer, pri izuchenii negamilp1tonova vzaimedee0stviya
(magnitnogo tipa) kvantovyh solitonov [22] ili vzaimodee0stviya
(super)strunnyh i strunnopodobnyh obp2ektov (naprimer, kvantovyh vihree0).
Odnako, korrektnaya postanovka pryamoe0 zadachi teorii predstavlenie0 v e1tom
sluchae avtoru neizvestna.

Vazhnoe0 nereshennoe0 problemoe0 yavlyaet\-sya globalizatsiya konstruktsii
izotopicheskih par, t.e. postroenie globalp1nyh algebraicheskih obp2ektov,
infinitezimalp1nye versii kotoryh predstavlyayut soboe0 izotopicheskie pary, a
takzhe razrabotka analoga teorii li dlya nih. Chastnye sluchai, svyazannye s
prostee0shchimi $r$--mat\-rich\-ny\-mi izotopicheskimi parami
($r$--mat\-rich\-ny\-mi magnitnymi algebrami Li), pozvolyayut predpolozhitp1
nalichie svyazi mezhdu ukazannym syuzhetom i teoriee0 kvantovyh grupp [23] i
kvantovyh $\tau$--funktsie0 [24].

Po-vidimomu, geometricheskie izotopicheskie pary mogut igratp1 opredelennuyu
rolp1 takzhe v formalizme asimptoticheskogo kvantovaniya [20].

Neobhodimo takzhe otmetitp1 vozmozhnye prilozheniya izotopicheskih par k
metodu obratnoe0 zadachi rasseyaniya v teorii nelinee0nyh
dif\-fe\-ren\-tsi\-alp1\-nyh uravnenie0 s chastnymi proizvodnymi [25] (sm.
takzhe v [26] iz\-lo\-zhe\-nie dlya sistem nelinee0nyh obyknovennyh
differentsialp1nyh uravnenie0) putem postroeniya (sr.[27]) assotsiirovannyh s
izotopicheskimi pa\-ra\-mi analogov izospektralp1nyh deformatsie0
(predstavlenie0 Laksa).

Izotopicheskie pary kak algebraicheskie0 apparat opisaniya sistem (kak
konechnyh, tak i beskonechnyh tsepochek) nekanonichski vzaimodee0stvuyushchih
obp2ektov [27] mogut bytp1 ispolp1zovany dlya opisaniya struktury i
vy\-yav\-le\-niya skrytyh simmetrie0 spektrov kvantovyh sistem (naprimer,
izotopicheskaya para nekanonicheski sparennyh ostsillyatorov [1,3,27] ---
spektrov ostsillyatornogo tipa [28]). Pri e1tom sushchestvennuyu rolp1 dolzhny
igratp1 obp2ekty bolee slozhnogo kombinatornogo tipa, chem opisannye vyshe
predstavleniya izotopicheskih par, takie kak, naprimer, {\cyi
graf--predstavleniya}, t.e. sovokupnosti otobrazhenie0
$\{T_1^{\alpha}, \alpha=1,\ldots N_1; T_1^{\beta}, \beta=1,\ldots N_2;
T_1^{\alpha}:V_1\mapsto\End(H)$, $T_2^{\beta}:V_2\mapsto\End(H)$ takih, chto
$$\aligned
T_1^{\alpha}([X,Y]_U)=&\sum_{\beta=1}^{N_2}\Cal P_{\alpha\beta}
(T_1^{\alpha}(X)T_2^{\beta}(U)T_1^{\alpha}(Y)-A_{XUY}
T_1^{\alpha}(Y)T_2^{\beta}(U)T_1^{\alpha}(X)),\\
T_2^{\beta}([U,V]_X)=&\sum_{\alpha=1}^{N_1}\Cal Q_{\alpha\beta}
(T_2^{\beta}(U)T_1^{\alpha}(X)T_2^{\beta}(V)-A_{UXV}
T_2^{\beta}(V)T_1^{\alpha}(X)T_1^{\beta}(U)),
\endaligned
$$
gde $X,Y\in V_1$, $U,V\in V_2$, $\Cal P$ i $\Cal Q$ --- dve matritsy
$N_1\times N_2$.

Nalichie podobnyh kombinatorno--netrivialp1nyh vysshih analogov\linebreak
predstavlenie0 yavlyaet\-sya otlichitelp1noe0 i spetsificheskoe0 po
sravneniyu s obychnymi algebrami chertoe0 algebraicheskih par.

\

\

\

\Refs
\roster
\item"[1]" {\cyie Yurp1ev D.V.}$/\!/$ {\cyre TMF. 1995. T.105. S.000-000}
[English electronic version (Texas Electronic Archive on Math. Phys.):
{\it mp\_arc/94-267\/} (1994)].
\item"[2]" {\cyie Berezin F.A.}, {\cyre Vvedenie v algebru i analiz ot
kommutiruyushchih i antikommutiruyushchih peremennyh. M., Nauka, 1983};
{\cyie Lee0tes D.A.}$/\!/$ {\cyre UMN. 1980. T.35, vyp.1. S.3-57};
{\cyie Manin Yu.I.}, {\cyre Kalibrovochnye polya i kompleksnaya geometriya.
M., Nauka, 1984}; {\cyie Ogievetskie0 V.I., Sokachev E.S.} / {\cyre Matem.
analiz 22. M., VINITI, 1984, S.137-174}; {\cyie Roslye0 A.A., Hudaverdyan O.M.,
Shvarts A.S.} / {\cyre Sovrem. probl. matem. Fundam. napravleniya 9, M.,
VINITI,
1986, S.247-284.}
\item"[3]" {\it Juriev D.}, Classical and quantum dynamics of noncanonically
coupled oscillators and Lie superalgebras / E--print (SISSA Electronic Archive
on Funct. Anal.): {\it funct-an/9409003} (1994), Russian J. Math. Phys.
[to appear]; On the dynamics of noncanonically coupled oscillators and its
hidden superstructure / Preprint ESI 167 (1994) and e--print (LANL Electronic
Archive on Solv. Integr. Systems): {\it solv-int/9503003} (1995).
\item"[4]" {\cyie Kurosh A.G.}, {\cyre Obshchaya algebra. M., 1974}.
\item"[5]" {\it Juriev D.}, On the nonHamiltonian interaction of two
rotators / E--print (MSRI Electronic Archive on Diff. Geom. and Global Anal.):
{\it dg-ga/9409004} (1994) and Report RCMPI/95-01 (1995).
\item"[6]" {\it Koecher M.}, Jordan algebras and their appliactions, 1962;
{\it Jacobson N.}, Structure and representations of Jordan algebras,
Providence,
1968; {\it Koecher M.}, An elementary approach to bounded symmetric domains,
Houston, 1969; {\cyie Loos O.}, {\cyre Simmetricheskie prostranstva, M., Nauka,
1985}.
\item"[7]" {\it Loos O.}, Jordan pairs. Springer--Verlag, 1975; {\cyie Kuzp1min
E.N., Shestakov I.P.\/}$/$ {\cyre Sov\-rem. probl. matem. Fundam. napravleniya
57. M., VINITI, 1992}.
\item"[8]" {\it Faulkner J.R., Ferrar J.C.\/}$/\!/$ Commun. Alg. 1980. V.8.
P.993-1013.
\item"[9]" {\it Juriev D.}, Topics in hidden symmetries / E--print (LANL
Electronic Archive on High Energy Phys.): {\it hep-th/9405050} (1994).
\item"[10]" {\cyie Yurp1ev D.V.}, {\cyre Harakteristiki par operatorov,
gibridy Li, skobki Puassona i nelinee0naya geometricheskaya algebra} /
E--print (SISSA Electronic Archive on Funct. Anal.): {\it funct-an/9411007}
(1994) and Report RCMPI/95-02 (1995).
\item"[11]" {\cyie Lee0tes D.A.} / {\cyre Sovrem. probl. matem. Novee0shie
dostizheniya 25, M., VINITI, 1984, S.3-50}.
\item"[12]" {\cyie Gisharde A.}, {\cyre Kogomologii topologicheskih grupp i
algebr Li. M., Mir, 1984}; {\cyie Fuks D.B.}, {\cyre Kogomologii
beskonechnomernyh algebr Li. M., Nauka, 1984}.
\item"[13]" {\cyie Zhelobenko D.P.\/}/\!/ {\cyre UMN. 1962. T.17, vyp.1.
S.27-120}; {\it Mickelsson J.\/}/\!/ Rep. Math. Phys. 1973. V.4. P.303-318,
1980. V.18. P.197-210; {\it Homberg A.\/}/\!/ Invent. Math. 1975. V.37.
P.42-47;
{\cyie Zhelobenko D.P.\/}/\!/ {\cyre Izv. AN SSSR, ser. matem. 1988. T.52.
S.758-773}; {\it Zhelobenko D.P.\/}/\!/ J. Group Theory in Phys. 1993. V.1.
P.201-233; {\cyie Zhelobenko D.P.}, {\cyre Predstavleniya reduktivnyh algebr
Li. M., Nauka, 1994}.
\item"[14]" {\cyie Asherova R.M., Smirnov Yu.F., Tolstoe0 V.N.\/}/\!/ {\cyre
TMF. 1971. T.8. S.255-271; Ma\-tem. zametki 1979. T.26. S.15-26}; {\cyie
Tolstoe0 V.N.\/}/\!/ {\cyre Teoretiko--gruppovye metody v fizike. 1986. T.2.
S.46-54}.
\item"[15]" {\it Okubo S.\/}$/\!/$ J. Math. Phys. 1994. V.35. P.2785-2803.
\item"[16]" {\cyie Semenov--Tyan--Shanskie0 M.A.\/}$/\!/$ {\cyre Funkts.
anal. i ego prilozh. 1983. T.13, vyp.4. S.17-33.}
\item"[17]" {\cyie Zhelobenko D.P., Shtern A.I.}, {\cyre Predstavleniya
grupp Li. M., Nauka, 1983}.
\item"[18]" {\cyie Bogolyubov N.N., Tolmachev V.V., Shirkov D.V.}, {\cyre
Novye0 metod v teorii sverhprovodimosti. Izd-vo AN SSSR, 1958.}
\item"[19]" {\cyie Barut A., Ronchka R.}, {\cyre Teoriya predstavlenie0 grupp
i ee prilozheniya. M., Mir, 1980.}
\item"[20]" {\cyie Karasev M.V., Maslov V.P.}, {\cyre Nelinee0nye skobki
Puassona. Geometriya i kvantovanie. M., Nauka, 1991.}
\item"[21]" {\it Nazaikinskii V., Sternin B., Shatalov V.}, Methods of
noncommutative analysis: theory and applications. Walter de Gruyter Inc., 1995.
\item"[22]" {\cyie Faddeev L.D., Korepin V.E.\/}$/\!/$ {\cyre TMF. 1975.
T.25. S.147-163}; {\it Faddeev L.D., Korepin V.E.\/}$/\!/$ Phys. Rep. C. 1978.
V.42, no.1.
\item"[23]" {\cyie Drinfelp1d V.G.\/}$/\!/$ {\cyre Zap. nauchn. semin. LOMI.
1986. T.155. S.18-49}; {\cyie Reshetihin N.Yu., Tahtadzhyan L.A., Faddeev
L.D.\/}
$/\!/$ {\cyre Algebra i anal. 1989. T.1. S.178-206}.
\item"[24]" {\it Gerasimov A., Khoroshkin S., Lebedev D., Morozov A.},
Generalized
Hirota equations and representation theory. I. The case of
$\operatorname{SL}(2)$
and $\operatorname{SL}_q(2)$ / E--print (LANL Electronic Archive on High Energy
Phys.): {\it hep-th/9405011\/} (1994); {\it Mironov A.}, Quantum deformations
of $\tau$--functions, bilinear identities and representation theory / E--print
(LANL Electronic Archive on High Energy Phys.): {\it hep-th/9409190\/} (1994);
{\it Kharchev S., Mironov A., Morozov A.}, Non--standard KP evolution and
quantum $\tau$--functions / E--print (Duke Univ. Electronic Archive on Q-Alg):
{\it q-alg/9501013\/} (1995).
\item"[25]" {\it Lax P.\/}$/\!/$ Commun. Pure Appl. Math. 1968. V.21.
P.467-490;
{\it Moser J.\/}$/\!/$ Adv. Math. 1975. V.16. P.197-220; {\cyie Zakharov V.E.,
Shabat A.B.\/}$/\!/$ {\cyre Funkts. anal. i ego prilozh. 1974. T.8, vyp.3.
S.43-53, 1979. T.13, vyp.3. S.13-22}; {\cyie Dubrovin B.A., Matveev V.B.,
Novikov S.P.\/}$/\!/$ {\cyre UMN. 1976. T.31, vyp.1. S.55-136;} {\cyie Zaharov
V.E., Manakov S.V., Novikov S.P., Pitaevskie0 L.P.}, {\cyre Teoriya solitonov.
Metod obratnoe0 zadachi. M., Nauka, 1980.}
\item"[26]" {\cyie Perelomov A.M.}, {\cyre Integriruemye sistemy klassicheskoe0
mehaniki i algebry Li. M., Nauka, 1990.}
\item"[27]" {\it Juriev D.\/}$/\!/$ Russian J. Math. Phys. 1995. V.3. no.4
[e--print version (LANL Electronic Archive on Solv. Integr. Systems):
{\it solv-int/9505001} (1995)].
\item"[28]" {\cyie Veselov A.P., Shabat A.B.\/}$/\!/$ {\cyre Funktsion.
anal. i ego prilozh. 1993. T.27, vyp.2. S.1-21}.
\endroster
\endRefs
\enddocument